\documentclass[a4paper]{amsart}

\newtheorem{thm}{Theorem}[section]

\newtheorem{prop}[thm]{Proposition}
\newtheorem{cor}[thm]{Corollary}

\theoremstyle{definition}
\newtheorem{defn}[thm]{Definition}

\theoremstyle{remark}

\begin{document}

\title[Galilean Connections]
{On Galilean connections and the first jet bundle}
\author{James D.E. Grant}
\author{Brad Lackey}
\address{Department of Mathematics,
University of Hull, Hull HU6 7RX, United Kingdom}
\email{J.D.Grant@maths.hull.ac.uk}
\email{B.Lackey@maths.hull.ac.uk}
\date{24 September 1999.}
\thanks{This work was funded by the EPSRC,
and the University of Hull.}

\begin{abstract}
We express the first jet bundle of curves in
Euclidean space as homogeneous spaces associated to a
Galilean-type group. Certain Cartan connections
on a manifold with values in the Lie algebra of the Galilean
group are characterized as geometries associated
to systems of second order ordinary differential equations.
We show these Cartan connections admit a form
of normal coordinates, and that in these normal coordinates
the geodesic equations of the connection are
second order ordinary differential equations. We then classify
such connections by some of their torsions,
extending a classical theorem of Chern involving the geometry
associated to a system of second order differential equations.
\end{abstract}
\maketitle

\section{A preparatory:
the Euclidean group and Riemannian geometry}

Recall that the Euclidean group of
motions can be described as a matrix group:
\begin{equation*}
Eucl_n = \left\{ \left(
\begin{array}{cc}
1 & 0 \\
x & A
\end{array} \right)\::\:
x \in \mathbb{R}^n, A \in
O_n\right\}\:.
\end{equation*}
The canonical left-invariant
Maurer-Cartan form for this group is
\begin{equation*}
\left(\begin{array}{cc}
1 & 0 \\
-A^{-1}x & A^{-1}
\end{array} \right)
\left(\begin{array}{cc}
0 & 0 \\
dx & dA
\end{array} \right)
= \left(
\begin{array}{cc}
0 & 0 \\
A^{-1}dx & A^{-1}dA
\end{array} \right)\:.
\end{equation*}
Writing this as
$\left(\begin{array}{cc}
0 & 0 \\
\omega^i & \omega^i{}_j
\end{array}\right)$,
one can easily see the Maurer-Cartan equations for the
component forms are the traditional structural equations of
Euclidean geometry:
\begin{eqnarray*}
d\omega^i + \omega^i{}_j \wedge \omega^j &=& 0\\
d\omega^i{}_j + \omega^i{}_k \wedge \omega^k{}_j &=& 0\:.
\end{eqnarray*}

Euclidean space is represented as a homogeneous space
associated to the Euclidean group,
$\pi:Eucl_n \to \mathbb{R}^n$. We can view
$Eucl_n$ as the total space of a principle
$O_n$-bundle over $\mathbb{R}^n$;
the Maurer-Cartan form is a
$\mathfrak{eucl}_n$-valued one form on
the total space of the bundle which satisfies
the Maurer-Cartan equations.

A Riemannian geometry is a manifold, $M$, and an
($O_n$-valued) orthonormal frame bundle, $\pi:P \to M$,
together with $\mathfrak{eucl}_n$-valued one form
on the total space of the bundle, generically called a
Cartan connection, but in this case
traditionally termed a metrical or
Riemannian connection. However, we do not require
(and it may not be possible) that this
matrix valued one-form satisfy the Maurer-Cartan equations.
This failure is measured as torsion/curvature of the forms:
\begin{equation*}
d\left(
\begin{array}{cc}
0 & 0 \\
\omega^i & \omega^i{}_j
\end{array}\right) +
\left(\begin{array}{cc}
0 & 0 \\
\omega^i & \omega^i{}_k
\end{array}\right) \wedge
\left(\begin{array}{cc}
0 & 0 \\
\omega^k & \omega^k{}_j
\end{array}\right) =
\left(\begin{array}{cc}
0 & 0 \\
\Omega^i & \Omega^i{}_j
\end{array}\right)\:.
\end{equation*}
We may write
\begin{eqnarray*}
\Omega^i &=& \frac{1}{2} T^i{}_{jk}\: \omega^j \wedge \omega^k\\
\Omega^i{}_j &=& \frac{1}{2} R^i{}_{jkl}\: \omega^k \wedge \omega^l\:.
\end{eqnarray*}
The tensors $T,R$ are the classical
Riemannian torsion and curvature.
The fundamental theorem of Riemannian
geometry can be phrased as: the tensor $T$
parameterizes in a one-to-one
fashion the metrical connections on
$P$. In particular, there exists a
unique connection such that $T = 0$.

\section{Klein and Cartan geometries}

The essence of Klein's Erlangen programme was to
realize the various forms of geometry that had arisen
in the nineteenth century as homogeneous spaces
associated to Lie groups. Specifically, these geometries
have transitive automorphism groups, $G$; so if
$H$ represents the isotropy group of a typical point,
the original space is diffeomorphic to $G/H$.
We have the canonical principle $H$-bundle $\pi:G \to
G/H$; the Maurer-Cartan forms of $G$, $\omega_G$,
is a $\mathfrak{g}$-valued one-form which satisfies
the Maurer-Cartan equations:
$d\omega_G + \frac{1}{2}[\omega_G,\omega_G] = 0$.

Cartan extended Klein's programme to include
curved versions of geometries. In modern language, we consider
an $H$-principle bundle over a smooth manifold,
$\pi:P \to M$, and a $\mathfrak{g}$-valued one-form, $\omega
\in \Omega^1(P,\mathfrak{g})$. This one-form,
called the Cartan connection, is to satisfy:
\begin{enumerate}
\item for each $p \in P$, $\omega|_p:TP \to \mathfrak{g}$ is an isomorphism;
\item for any $h \in H$, $R^*_h\omega = ad(h^{-1})\:\omega$; and,
\item for each $X \in \mathfrak{h}$, $\omega(X^\dagger) = X$.
\end{enumerate}
These three conditions guarantee that at each
$x \in M$, $T_xM \cong \mathfrak{g}/\mathfrak{h}$.
Since the Maurer-Cartan equations are generally not satisfied,
$M$ can be viewed as a curved version of $G/H$.

Fix $\pi:P \to M$ and a Cartan connection $\omega$.
For any curve $\rho:[0,1] \to P$, the equation
$\rho^*\omega = \tilde{\rho}^*\omega_G$ is a
first order ordinary differential equation for the curve
$\tilde{\rho}$ into $G$. Specifying an initial condition,
say $\tilde{\rho}(0) = e$, and using the
compactness of $[0,1]$, one gets a unique curve
$\tilde{\rho}:[0,1] \to G$ called the
\textit{development} of $\rho$. Now, given a curve
$\sigma:[0,1] \to M$, let $\rho$ be any lift of
$\sigma$ to $P$. One can show that the development,
$\tilde{\rho}$, depends upon the lift, but the
projection of this curve to $G/H$ does not.
Therefore, one gets a unique curve
$\tilde{\sigma}:[0,1] \to G/H$ associated to $\sigma$,
also called its development. The purpose of
this exercise is that if one has a natural notion
of `straight line' in the model space $G/H$,
then one can define geodesics in $M$ by taking
curves whose developments are straight lines.
This is precisely the way geodesics arise in
Riemannian and projective geometry.
See \cite[\S\S{\textbf{3}.7 and \textbf{5}.4}]{S}
for a more detailed treatment of developments of curves.

\section{The first jet bundle}

Our goal is to realize the first jet bundle
$J^1(\mathbb{R},\mathbb{R}^n)$ as a homogeneous space
associated to its natural automorphism group.
To this end, we introduce some features of the first jet
bundle.

Formally, the first jet bundle, $J^1(\mathbb{R},\mathbb{R}^n)$,
is a vector bundle over $\mathbb{R} \times
\mathbb{R}^n$ defined via equivalence classes
of germs of curves from $\mathbb{R}$ into $\mathbb{R}^n$.
Although it is much easier to understand the
structure of the \textit{first} jet bundle by using local
coordinates, for completeness we present the
formal definition as this is simple to generalize.

\begin{defn}
Fix $(t,x) \in \mathbb{R} \times \mathbb{R}^n$ and consider two
germs of curves $\gamma_1,\gamma_2: (t-\epsilon,t+\epsilon) \to
\mathbb{R}^n$ with $\gamma_1(t) = x = \gamma_2(t)$. We say these
curves have first order contact at $(t,x)$ if they have the same
tangent at $(t,x)$: $\dot{\gamma}_1(t) = \dot{\gamma}_2(t)$. First
order contact is clearly an equivalence relation, so denote
$[\gamma]_{(t,x)}$ as the equivalence class of $\gamma$. The first
jet bundle $J^1(\mathbb{R},\mathbb{R}^n)$ is the vector bundle
over $\mathbb{R} \times \mathbb{R}^n$ whose fibre over $(t,x)$ is
$\{[\gamma]_{(t,x)}\: |\: \gamma: (t-\epsilon,t+\epsilon) \to
\mathbb{R}^n, \ \gamma(t) = x\}$.
\end{defn}

It is straight forward to justify that the definition above does
indeed yield $J^1(\mathbb{R},\mathbb{R}^n)$ as a vector bundle
over $\mathbb{R} \times \mathbb{R}^n$. One way to see this is to
identify $J^1(\mathbb{R},\mathbb{R}^n) \cong \mathbb{R} \times
T\mathbb{R}^n$. At each $(t,x)$, an equivalence class
$[\gamma]_{(t,x)}$ determines a unique tangent vector to
$\mathbb{R}^n$. Conversely, each tangent vector to $\mathbb{R}^n$
at $(t,x)$ can be realized as the tangent to a curve through
$(t,x)$.

The identification $J^1(\mathbb{R},\mathbb{R}^n) \cong \mathbb{R}
\times T\mathbb{R}^n$ yields a preferred choice of coordinates.
Specifically, if $(x^j)$ is a local coordinate system of
$\mathbb{R}^n$, then tangent vectors are expressed as $y =
y^j\frac{\partial}{\partial x^j}$. So we may always take local
coordinates of the form $(t,x^j,y^j)$. Note that if we change
coordinates, $\bar{x} = \bar{x}(x)$, we have both $d\bar{x} =
A^{-1}\:dx$ and $\bar{y} = A^{-1}\:y$, where $A = \frac{\partial
x}{\partial \bar{x}}$. Thus the `contact' forms $dx - y\:dt$ have
nice transformational character.

One cannot overestimate the utility of the contact forms. If we
represent a curve in coordinates, $t \stackrel{\gamma}{\mapsto}
x(t)$, then we can naturally define a curve on
$J^1(\mathbb{R},\mathbb{R}^n)$ by $t \mapsto
(t,x(t),\frac{dx}{dt}(t))$. This curve is called the first jet (or
prolongation) of $\gamma$, and is traditionally denoted
$j^1\gamma$. Now, if $\sigma:\mathbb{R} \to
J^1(\mathbb{R},\mathbb{R}^n)$, then $\sigma = j^1 \gamma$ for some
$\gamma$ if and only if $\sigma$ is integral to the contact forms:
$\sigma^*(dx - y\:dt) = 0$. Thus the contact forms completely
encode the fact that $y$ is the derivative of $x$ with respect to
$t$.

The first jet bundle, like Euclidean space or projective space,
carries a canonical notion of straight line. Intuitively, a
straight line is given by $x = at$ and $y = a$ for some fixed $a
\in \mathbb{R}^n$. Of course, if $a=0$, then the curve degenerates
to just a point, and we wish to exclude this. We can express this
picture using the contact forms together with $dt,dy$.

\begin{defn}
A curve $\tilde{\sigma}:\mathbb{R} \to
J^1(\mathbb{R},\mathbb{R}^n)$ is a \em{straight line} if
$\tilde{\sigma}^*(dx - y\:dt) = \tilde{\sigma}^*(dy) = 0$ and
$\tilde{\sigma}^*(dt) \not= 0$.
\end{defn}

\section{Galilean geometry}

Strictly speaking we are not considering the Galilean group;
rather than taking the rotations, translations, and boosts on the
tangent space, we take arbitrary affine transformations and boosts
on the first jet space.\footnote{In what follows, one can use the
traditional Galilean group, but the resulting geometry is slightly
different. It is, in fact, that proposed in \cite{L2}.}
Specifically, we define
\begin{equation*}
Gal_n = \left\{\left(
\begin{array}{ccc}
1 & 0 & 0 \\
t & 1 & 0 \\
x & y & A
\end{array} \right)\::\:
t \in \mathbb{R},
x,y \in \mathbb{R}^n,
A \in Gl_n\right\}\:.
\end{equation*}
The canonical Maurer-Cartan forms for this group are
\begin{eqnarray*}
&&\left(
\begin{array}{ccc}
1 & 0 & 0 \\
-t & 1 & 0 \\
-A^{-1}(x - yt) &
-A^{-1}y & A^{-1}
\end{array}
\right) \left(
\begin{array}{ccc}
0 & 0 & 0 \\
dt & 0 & 0 \\
dx & dy & dA
\end{array}\right)\\ && \qquad
\qquad = \left(
\begin{array}{ccc}
0 & 0 & 0 \\
dt & 0 & 0 \\
A^{-1}(dx - y\:dt) & A^{-1}dy & A^{-1}dA
\end{array}\right)\:.
\end{eqnarray*}

The isotropy group of a point in the first jet bundle is
\begin{equation*}
H = \left\{ \left(
\begin{array}{ccc}
1 & 0 & 0 \\
0 & 1 & 0 \\
0 & 0 & A
\end{array}\right)\::\:
A \in Gl_n\right\}\:.
\end{equation*}
Clearly $J^1(\mathbb{R},\mathbb{R}^n) \cong Gal_n/H$.

\begin{defn}
A Galilean manifold is a smooth manifold $X$,
equipped with an principle $H$-bundle, $\pi:P \to X$, and a
Cartan connection on $P$ with values in $\mathfrak{gal}_n$.
\end{defn}

To fix our notation, we will write the
Cartan connection of our Galilean manifold as
$\left(
\begin{array}{ccc}
0 & 0 & 0 \\
\tau & 0 & 0 \\
\omega & \phi & \Pi
\end{array}\right)$. The
curvature of our connection is then
\begin{equation*}
\left(
\begin{array}{ccc}
0 & 0 & 0 \\
T & 0 & 0 \\
\Omega & \Phi & R
\end{array}\right) =
\left(
\begin{array}{ccc}
0 & 0 & 0 \\
d\tau & 0 & 0 \\
d\omega + \Pi \wedge \omega + \phi \wedge \tau &
d\phi + \Pi \wedge \phi & d\Pi + \Pi \wedge \Pi
\end{array}\right)\:.
\end{equation*}
Curvature measures the failure of the geometric structure to be
representable in coordinate systems. In the first jet bundle, we
may always take adapted coordinates of the form $(t,x,y)$; but for
the most general Galilean geometries this is not possible even
locally. Therefore, we are really interested in Cartan connections
which are ``partially integrable.''

\begin{defn}
The Cartan connection of a Galilean manifold, $X$, is
\em{holonomic} if $T=0$ and \em{torsion-free} if $\Omega = 0$. If
both conditions are satisfied, we will call $X$ a normal Galilean
manifold.
\end{defn}

We use the overworked term `normal' for this type of Galilean
manifold because it has normal coordinates in some sense. The
following theorem is a precise statement of this; intuitively, we
have demanded certain curvatures vanish in order to find local
coordinates which mimic the adapted coordinates of the first jet
bundle.

\begin{thm}
Let $X$ be a normal Galilean manifold. Take any point in $X$. Then
about this point, there exists a coordinate system, $(t,x,y)$, and
a local section $h$ of $P$, such that $h^*\tau = dt$, $h^*\omega =
dx - y\:dt$, and $h^*\phi = dy$ (mod $dt,dx$).
\end{thm}
\begin{proof}
Let $\tilde{h}$ be any local section of $P$ near our given point.
Immediately we have $d\tilde{h}^*\tau = 0$, and hence restricting
our domain if necessary $\tilde{h}^*\tau = dt$. As $dt \not=0$,
the function $t$ serves as a coordinate function.

Define $\mathcal{I}$ to be the exterior ideal generated by the
one-forms $dt$ and $\tilde{h}^*\omega$. The curvature conditions
$T = 0$ and $\Omega = 0$ show that $\mathcal{I}$ is integrable.
Since $dt \in \mathcal{I}$ each integral manifold lies in a
$t=\mathrm{constant}$ hypersurface. Again restricting our domain
if necessary, we may choose independent functions $x^j$ so that
the integral manifolds of $\mathcal{I}$ are uniquely described by
$t=\mathrm{constant}$ and $x =\mathrm{constant}$.

Note that the integral manifolds of $\mathcal{I}$ are also the
integral manifolds of $\langle dt, dx^j \rangle$. Therefore, at
each point of our domain we may express $\tilde{h}^*\omega = A(dx
- y\:dt)$ for some invertible matrix $A$ and vector $y$. We claim
that $\{y^j\}$ serve as coordinates of the integral manifolds, and
hence completes our local coordinate system. To show this, it
suffices to show that $\{dt,dx,dy\}$ forms a local coframe near
our given point. To this end, we compute:
\begin{equation*}
d\tilde{h}^*\omega =
dA\cdot A^{-1} \wedge \tilde{h}^*\omega - A\:dy \wedge dt\:.
\end{equation*}
Yet, as $\Omega = 0$, we have the structure equation $d\omega =
-\Pi \wedge \omega - \phi \wedge \tau$. Therefore,
\begin{equation*}
(dA\cdot A^{-1} + \tilde{h}^*\Pi) \wedge \tilde{h}^*\omega =
(A\:dy - \tilde{h}^*\phi) \wedge dt\:.
\end{equation*}
Every two-form in the left side of this equation contains an
$\tilde{h}^*\omega$. Thus, $A\:dy - \tilde{h}^*\phi =
N\tilde{h}^*\omega + \lambda\:dt$ for some matrix $N$ and vector
$\lambda$. Yet, $\tau, \omega,\phi$ are linearly independent.
Therefore, the $dy^j$ are linearly independent, and also
independent from $\tilde{h}^*\omega, \tilde{h}^*\tau$, and thus
independent from $dt,dx$ as well.

Finally, define a new section by $h = R_A\tilde{h}$, where $A$
abbreviates $\left(\begin{array}{ccc} 1 & 0 & 0 \\ 0 & 1 & 0 \\ 0
& 0 & A \end{array}\right) \in H$. Then we have $h^*\tau =
\tilde{h}^*\tau$, $h^*\omega = A^{-1}\tilde{h}^*\omega$, and
$h^*\phi = A^{-1}\tilde{h}^*\phi$. In particular, $h^*\tau = dt$,
$h^*\omega = dx - y\: dt$, and $h^*\phi = dy$ (mod $h^*\tau,
h^*\omega$), as desired.
\end{proof}

\begin{defn}
A \em{geodesic} in a Galilean manifold is a (germ of a) curve
$\sigma:(-a,a) \to X$ whose development $\tilde{\sigma}:(-a,a) \to
J^1(\mathbb{R},\mathbb{R}^n)$ is contained in a straight line.
\end{defn}

\begin{prop}
Let $X$ be a normal Galilean manifold, and let $(t,x,y)$ be normal
coordinates in $X$, and write $h^*\phi = dy + \Gamma(t,x,y)\:dt +
N(t,x,y)(dx - y\: dt)$. Then the geodesics of $X$ are the
solutions to the system of second order differential equations
\begin{equation}
\frac{d^2x^j}{dt^2} + \Gamma^j\left(t,x,\frac{dx}{dt}\right) = 0\:.
\label{ode}
\end{equation}
\end{prop}
\begin{proof}
We have already seen that the straight lines in
$J^1(\mathbb{R},\mathbb{R}^n)$ are characterized by the equations
$\tilde{\sigma}^*(dx - y\: dt) = 0$ and $\tilde{\sigma}^*(dy) =
0$. Yet, $dx - y\:dt$ and $dy$ are the $\omega$ and $\phi$ terms,
respectively, in the Maurer-Cartan forms of the Galilean group.
Thus by the definition of development, a curve $\sigma$ is a
geodesic if and only if it is non-degenerate and $\sigma^* (h^*
\omega) = \sigma^* (h^* \phi) = 0$. Yet, in normal coordinates
$h^*\omega = dx - y\:dt$ and $h^*\phi = dy + \Gamma\:dt + N(dx -
y\: dt)$; the integral curves to these forms are precisely the
solutions to (\ref{ode}).
\end{proof}

\begin{thm}
Consider local coordinates $(t,x,y)$, and curves in this
coordinate system given by (\ref{ode}). Then we have the following
\begin{enumerate}
\item[(i)] There are many normal Galilean geometries
having these coordinates as normal coordinates and
(\ref{ode}) as geodesic equations. For these we have
\begin{equation*}
\Phi^i = D^i{}_j \tau \wedge \phi^j + Q^i{}_{jk} \omega^j \wedge
\phi^k + P^i{}_j \tau \wedge \omega^j + \frac{1}{2}T^i{}_{jk}
\omega^j \wedge \omega^k\:.
\end{equation*}
\item[(ii)] For the geometries described in (i),
the tensor $D$ parameterizes in a one-to-one fashion
the forms $\phi$. There are many such normal Galilean geometries
for each choice of $D$.
\item[(iii)] To each geometry described in (ii),
the symmetric part of the tensor $Q$,
$Q^i{}_{(jk)} = \frac{1}{2}(Q^i{}_{jk} + Q^i{}_{kj})$,
parameterizes in a one-to-one fashion the forms $\Pi$. There is
precisely one such normal Galilean geometry for each choice of
$Q^i{}_{(jk)}$.
\end{enumerate}
\end{thm}
\begin{proof}
Let $(t,x,y)$ be as given; on the portion of $P$ over this
coordinate system we have natural coordinates $(t,x,y,A)$. {}From
the construction of the coordinate $t$, it is clear that $\tau =
dt$. Similarly, if there is a section $h$ so that
$h^*\omega = dx - y\:dt$, then we must have
$\omega = A^{-1}(dx - y\:dt)$
(explicitly, the desired section is then $h(t,x,y) =
(t,x,y,A(t,x,y))$). This implies that $\phi = A^{-1}dy$ (mod
$\tau, \omega$). The geodesics of the geometry are given by
$\sigma^*\omega = \sigma^*\phi = 0$, hence $\phi = A^{-1}(dy +
\Gamma\:dt)$ (mod $\omega$), say $\phi = A^{-1}(dy + \Gamma\:dt+
N(dx - y\:dt))$ where $N$ is an arbitrary matrix. If the structure
equation $d\omega = -\Pi \wedge \omega - \phi \wedge \tau$ is to
be satisfied, we must have $\Pi^i{}_j = (A^{-1}dA)^i{}_j +
(A^{-1})^i{}_{k}N^k{}_{l}A^l{}_j\:dt + (A^{-1})^i_{\
l}\Gamma^l{}_{rs}A^r{}_jA^s{}_{k}\:\omega^k$, where $\Gamma^i_{\
jk} = \Gamma^i{}_{kj}$ but is otherwise arbitrary. Hence (i) is
shown. Moreover, the most general Galilean geometry satisfying (i)
is
\begin{eqnarray*}
\tau &=& dt\\ \omega^i &=& (A^{-1})^i{}_j(dx^j - y^j\:dt)\\
\phi^i &=& (A^{-1})^i{}_j (dy^j + \Gamma^j\:dt + N^j{}_{k}(dx^k
- y^k\:dt))\\ \Pi^i{}_j &=& (A^{-1}dA)^i{}_j + (A^{-1})^i_{\
k} N^k{}_{l}A^l{}_j\:dt + (A^{-1})^i{}_{l}\Gamma^l{}_{rs}A^r_{\
j}A^s{}_{k}\:\omega^k\:.
\end{eqnarray*}
We leave it as an easy exercise to the reader to show that for
these forms, $\Phi$ has no $\phi \wedge \phi$ term.

To show (ii), we compute
\begin{eqnarray*}
\Phi^i &=& d\phi + \Pi \wedge \phi \\ &=& (A^{-1})^i{}_j
\left(
2N^j{}_{k} - \frac{\partial \Gamma^j}{\partial y^k}
\right)
A^k{}_l \: \tau \wedge \phi^l \quad (\text{mod } \omega)\:.
\end{eqnarray*}
Therefore, for any choice of $D$, we have
$N = \frac{1}{2}\left(
\frac{\partial \Gamma}{\partial y} +
ADA^{-1}\right)$.

Finally, suppose $D$, and hence $N$, is now fixed.
Note that we may write\newline
$\phi = A^{-1}(dy + \Gamma\:dt
+ \frac{1}{2}\frac{\partial \Gamma}{\partial y} (dx - y\: dt))
+ D\:\omega$. One easily computes the
$\omega^j \wedge \phi^k$ term in $\Phi^i$ to be
\begin{equation*}
(A^{-1})^i{}_{l}\left(\Gamma^l{}_{rs} - \frac{1}{2}
\frac{\partial^2 \Gamma^l}{\partial y^r \partial y^s} -
\frac{1}{2} A^l{}_{m} \frac{\partial D^m{}_{p}}{\partial y^s}
(A^{-1})^p{}_{r} \right)A^r{}_jA^s{}_{k}\: \omega^j \wedge
\phi^k\:.
\end{equation*}
Only the symmetric part of this involves $\Gamma^i{}_{jk}$, hence
specifying $Q^i{}_{(jk)}$ determines $\Gamma^i{}_{jk}$ uniquely,
and (iii) is proven.
\end{proof}

This theorem generalizes the classical theorem of Chern about the
geometry induced by systems of second order ordinary differential
equation. We can state his result in our language as the
following.

\begin{cor}[Chern \cite{C}]
To each system of second order ordinary differential equations,
(\ref{ode}), there exists a unique normal Galilean connection such
that $\Phi = 0$ (mod $\omega$).
\end{cor}
\begin{proof}
Choosing $D^i{}_j = 0$ and $Q^i{}_{(jk)} = 0$,
we see that $\Phi = 0$ (mod $\omega$).
\end{proof}

\appendix
\section{Geometry of jet bundles}

Jet spaces in generality have a great deal of structure. In recent
years, it has become popular to study the geometry of jet bundles,
see for instance \cite{Sa}; however, the usual techniques involve
nonlinear connections in the tangent bundle of the jet bundle.
Throughout the paper we have used the exterior differential
calculus, as this formalism is very convenient for studying Cartan
geometries. In this appendix, we indicate how the results of the
paper impact upon the geometry of jet spaces.

Let $B$ be a smooth manifold, and suppose $\pi:X \to B$ is a fibre
bundle with fibre $F$. The tangent vectors to the fibres have good
transformational character, so form a vector subbundle of $TX$
called the vertical vector space, and denoted as $VTX$. The easy
way to see this is to note that $\xi$ is tangent to the fibre if
and only if $\pi_*\xi = 0$, hence $VTX = \mathrm{ker} \pi_*$.
Unfortunately, although vectors tangent to $B$ locally complement
the vertical vector space, they do not transform properly on $X$.
Thus $TB$ is not a subbundle of $TX$. A \textit{nonlinear
connection} is merely a selection of complementary vector bundle
to $VTX$ in $TX$. This bundle is usually called the horizontal
vector space, and denoted $HTX$.

Another definition of a nonlinear connection is via the short
exact sequence of vector bundles over $X$:
\begin{equation*}
0 \to VTX \to TX \stackrel{\pi_*}{\to} \pi^*TB \to 0\:,
\end{equation*}
where $\pi^*TB$ is the pull-back of the tangent of $B$ over $X$. A
nonlinear connection is a splitting of this short exact sequence:
$TX \cong VTX \oplus HTX$ where $HTX \cong \pi^*TB$.
Unfortunately, without further knowledge about the structure of
$X$, there is no canonical splitting.

Now we specialize to the case where $B = \mathbb{R} \times M$ for
some smooth manifold $M$, and $X = J^1(\mathbb{R},M)$. We wish to
describe a nonlinear connection in local adapted coordinates
$(t,x,y)$. The vertical vector space over this coordinate chart is
spanned by $\left\{\frac{\partial}{\partial y^j}\right\}$, and
$\pi^*T(\mathbb{R}\times M)$ is spanned by
$\left\{\frac{\partial}{\partial t} + y^j\frac{\partial}{\partial
x^j}, \frac{\partial}{\partial x^j}\right\}$. Do note that our
choice of vector fields in the latter contains the contact
structure of $J^1(\mathbb{R},M)$.

Over our coordinate system, a selection of nonlinear connection is
just a selection of functions $\Gamma,N$ such that
\begin{equation*}
HTJ^1(\mathbb{R},M) = \mathrm{span} \left\{
\frac{\partial}{\partial t} + y^j\frac{\partial}{\partial x^j} -
\Gamma^j(t,x,y)\frac{\partial}{\partial y^j},
\frac{\partial}{\partial x^j} - N^j_{\
k}(t,x,y)\frac{\partial}{\partial y^k} \right\}\:.
\end{equation*}
For convenience, let us adopt the notation $\frac{d}{dt} =
\frac{\partial}{\partial t} + y^j\frac{\partial}{\partial x^j} -
\Gamma^j(t,x,y)\frac{\partial}{\partial y^j}$ and
$\frac{\delta}{\delta x^j} = \frac{\partial}{\partial x^j} -
N^j{}_{k}(t,x,y)\frac{\partial}{\partial y^k}$.

We see that the nonlinear connection is expressed as two sets of
function. The functions $\Gamma$ describe a second order ordinary
differential equation, and as we have seen, the geometry of the
first jet bundle is closely tied to these differential equations.
One goal of studying this type of geometry is to determine a
natural choice of the functions $N$ based upon a selection of
$\Gamma$. The most natural thing to do is compute the commutators
of the vector fields above. Oddly enough, this is not the proper
thing to do.

Instead, it is far better to study affine connections on
$TJ^1(\mathbb{R},M)$. We are not interested in arbitrary affine
connections, but rather ones which reflect the structure of the
first jet space. In other words we are not interested in Cartan
connections with values in $\mathfrak{aff}_{2n+1}$ but rather
those with values in $\mathfrak{gal}_n$. That is, a normal
Galilean connections yields a covariant derivative with structure:
\begin{equation*}
\begin{array}{rclcrclcrcl}
\nabla_{\frac{d}{dt}} \frac{d}{dt} &=& 0 &\quad&
\nabla_{\frac{\delta}{\delta x^k}} \frac{d}{dt} &=& 0 &\quad&
\nabla_{\frac{\partial}{\partial y^k}} \frac{d}{dt} &=& 0 \\
\nabla_{\frac{d}{dt}} \frac{\delta}{\delta x^k} &=& N^j{}_{k}
\frac{\delta}{\delta x^j} &\quad& \nabla_{\frac{\delta}{\delta
x^k}} \frac{\delta}{\delta x^k} &=& \Gamma^l{}_{jk}
\frac{\delta}{\delta x^l} &\quad& \nabla_{\frac{\partial}{\partial
y^k}} \frac{\delta}{\delta x^k} &=& 0 \\ \nabla_{\frac{d}{dt}}
\frac{\partial}{\partial y^k} &=& N^j{}_{k}
\frac{\partial}{\partial y^j} &\quad& \nabla_{\frac{\delta}{\delta
x^k}} \frac{\partial}{\partial y^k} &=& 0 &\quad&
\nabla_{\frac{\partial}{\partial y^k}} \frac{\partial}{\partial
y^k} &=& \Gamma^l{}_{jk} \frac{\partial}{\partial y^k} \:.
\end{array}
\end{equation*}

Now, it is a simple matter to see that the components of torsion
determine the nonlinear and affine connection coefficients. For
instance,
\begin{equation*}
\nabla_{\frac{\partial}{\partial y^k}} \frac{d}{dt} -
\nabla_{\frac{d}{dt}} \frac{\partial}{\partial y^k} -
\left[\frac{\partial}{\partial y^k}, \frac{d}{dt}\right] = -
\frac{\delta}{\delta x^k} + \left(\frac{\partial
\Gamma^j}{\partial y^k} - 2N^j{}_{k}\right)
\frac{\partial}{\partial y^j}\:.
\end{equation*}

\end{document}